\documentclass[12pt]{amsart}
\usepackage[matrix,arrow]{xy}
\usepackage{amssymb,amsfonts,amsmath,footnote}

\newtheorem{theorem}{Theorem}[section]
\newtheorem{proposition}[theorem]{Proposition}
\newtheorem{lemma}[theorem]{Lemma}

\newtheorem{corollary}[theorem]{Corollary}
\newtheorem{remark}[theorem]{Remark}

\begin{document}

\def\colim{\mathop{\mathrm{colim}}}
\def\la#1{\mathop{\longleftarrow}\limits^{#1}}
\def\ra#1{\mathop{\longrightarrow}\limits^{#1}}
\def\be{\begin{equation}}
\def\ee{\end{equation}}

\def\VV{\mathbb{V}}
\def\MMU{\mathbb{MU}}
\def\MMO{\mathbb{MO}}
\def\HH{\mathbb{H}}
\def\SS{\mathbb{S}}
\def\LL{\mathbb{L}}
\def\FF{\mathbb{F}}
\def\RR{\mathbb{R}}
\def\ZZ{\mathbb{Z}}
\def\CC{\mathbb{C}}
\def\NN{\mathbb{N}}

\def\uZZ{\underline{\mathbb{Z}}}

\def\ev{\mathrm{ev}}
\def\map{\mathrm{map}}
\def\Hom{\mathrm{Hom}}
\def\End{\mathrm{End}}
\def\Ext{\mathrm{Ext}}
\def\Sq{\mathrm{Sq}}
\def\deg{\mathrm{deg}}
\def\id{\mathrm{id}}
\def\Th{\mathrm{\mu}}
\def\ker{\mathrm{ker}\,}
\def\Ann{\mathrm{Ann}}
\def\Mon{\mathrm{Mon}}
\def\res{\mathrm{Res}}
\def\fr{\mathrm{fr}}
\def\Sk{\mathrm{Sk}}
\def\Aut{\mathrm{Aut}}
\def\ob{\mathrm{ob}\,}
\def\Fun{\mathrm{Fun}}
\def\Stab{\mathrm{Stab}}
\def\Ind{\mathrm{Ind}}
\def\Sym{\mathrm{Sym}}
\def\op{\mathrm{op}}

\def\tensor{\otimes}

\def\cS{\mathcal{S}}
\def\cA{\mathcal{A}}
\def\cM{\mathcal{M}}
\def\cI{\mathcal{I}}
\def\cT{\mathcal{T}}
\def\cC{\mathcal{C}}
\def\cF{\mathcal{F}}
\def\cP{\mathcal{P}}

\def\hS{\widehat{S}}
\def\hW{\widehat{W}}

\def\MR{\MMU_\RR}
\def\MRn{\MMU^{(n)}}

\def\ovr{\smash{\overline{r}}}
\def\orho{\smash{\overline{\rho}}}
\def\ofd{\smash{\overline{\mathfrak{d}}}}
\def\oN{\smash{\overline{N}}}
\def\cA{\mathcal{A}}
\def\Sq{\mathrm{Sq}}
\def\Span{\mathrm{Span}}

\date{\today}

\title[The anti-automorphism of the Steenrod algebra]
{On the anti-automorphism \\
of the Steenrod algebra: II}

\author{V.~Giambalvo}
\address{Department of Mathematics\\
University of Connecticut\\
Storrs, CT 06269}
\email{vince@math.uconn.edu}

\author{H.~R.~Miller}
\address{Department of Mathematics\\
Massachusetts Institute of Technology\\
Cambridge, MA 02139}
\email{hrm@math.mit.edu}



\begin{abstract}
The relations of Barratt and Miller are shown to include all relations among 
the elements $P^i\chi P^{n-i}$ in the mod $p$ Steenrod algebra, and a minimal
set of relations is given.
\end{abstract}

\maketitle

\section{Introduction}

Milnor \cite{milnor} 
observed that the mod 2 Steenrod algebra $\cA$ forms a Hopf algebra with 
commutative diagonal determined by
\be
\Delta \Sq^n=\sum_i\Sq^i\tensor\Sq^{n-i}\,.
\ee
This allowed him to interpret the Cartan formula
as the assertion that the cohomology of a space forms a module-algebra
over $\cA$. The anti-automorphism $\chi$ in the Hopf algebra structure,
defined inductively by
\be
\chi\Sq^0=\Sq^0\quad,\quad
\sum_i\Sq^i\chi\Sq^{n-i}=0\quad\hbox{for}\quad n>0\,,
\label{def-chi}
\ee
has a topological interpretation too: If $K$ is a finite complex then
the homology of the Spanier-Whitehead dual $DK_+$ of $K_+$ is canonically
isomorphic to the cohomology of $K$. Under this isomorphism the left
action by $\theta\in\cA$ on $H^*(K)$ corresponds to the right action of 
$\chi\theta\in\cA$ on $H_*(DK_+)$.

In 1974 Davis \cite{davis}
proved that sometimes much more efficient ways exist to
compute $\chi\Sq^n$:
\be
\chi\Sq^{2^r-1}=\Sq^{2^{r-1}}\chi\Sq^{2^{r-1}-1}
\label{davis-one}
\ee
\be
\chi\Sq^{2^r-r-1}=\Sq^{2^{r-1}-1}\chi\Sq^{2^{r-1}-r}
+\Sq^{2^{r-1}}\chi\Sq^{2^{r-1}-r-1}
\label{davis-two}
\ee
Similarly, Straffin \cite{straffin} proved that if $r\geq0$ and $b\geq2$ then
\be
\sum_i\Sq^{2^ri}\chi\Sq^{2^r(b-i)}=0\,.
\ee
Both authors give analogous identities among reduced powers and their 
images under $\chi$ at an odd prime as well.

Barratt and Miller \cite{barratt-miller} found a general family of identities
which includes all these as special cases. We state it for the general prime.
When $p=2$, $P^n$ denotes $\Sq^n$. Let $\alpha(n)$ denote the sum of the 
$p$-adic digits of $n$. 

\begin{theorem} \cite{barratt-miller,crabb-crossley-hubbuck}
For any integer $k$ and any integer $l\geq0$ such that $pl-\alpha(l)<(p-1)n$
\be
\sum_i\binom{k-i}{l}P^i\chi P^{n-i}=0\,.
\label{bm}
\ee
\label{th-bm}
\end{theorem}

The defining relations occur with $l=0$. 
Davis's formulas (for $p=2$) are the cases in which 
$(n,l,k)=(2^r-1,2^{r-1}-1,2^r-1)$ or $(n,l,k)=(2^r-r-1,2^{r-1}-2,2^r-2)$.
Straffin's identities (for $p=2$) occur as $(n,l,k)=(2^rb,2^r-1,2^rb)$.

Since $\binom{(k+1)-i}{l}-\binom{k-i}{l}=\binom{k-i}{l-1}$,
the cases $(l,k+1)$ and $(l,k)$ of  (\ref{bm}) imply it for
$(l-1,k)$. Thus the relations for 
\be
l=l(n)=\max\{j:pj-\alpha(j)<(p-1)n\}
\label{ln}
\ee
imply all the rest. When $p=2$,
$l(2^r-1)=2^{r-1}-1$ and $l(2^r-r-1)=2^{r-1}-2$,
so Davis's relations are among these basic relations. 

Two questions now arise. To express them uniformly in the prime, let 
$\cP$ denote the algebra of Steenrod reduced powers (which is the full
Steenrod algebra when $p=2$), but assign $P^n$ degree $n$. It is natural to
ask:  \\
-- Are there yet other linear relations among the $n+1$
elements $P^i\chi P^{n-i}$ in $\cP^{n}$? \\
-- What is a minimal spanning set for 
\[
V_n=\Span\{P^i\chi P^{n-i}:0\leq i\leq n\}\subseteq\cP^n
\]
We answer these questions in
Theorem \ref{main-theorem} below.

Write $e_i,0\leq i\leq n$, for the $i$th standard basis vector in 
$\FF_p^{n+1}$.

\begin{proposition} 
For any integers $l,m,n$, with $0\leq l\leq n$,  
\be
\left\{\sum_i\binom{k-i}{l}e_i:m\leq k\leq m+l\right\}
\label{basis}
\ee
is linear independent in $\FF_p^{n+1}$.
\label{prop-basis}
\end{proposition}

\begin{proposition}
The set 
\be
\left\{P^i\chi P^{n-i}:l(n)+1\leq i\leq n\right\}
\label{basis2}
\ee
is linearly independent in $\cP^n$.
\label{prop-basis2}
\end{proposition}

Define a linear map
\be
\mu:\FF_p^{n+1}\rightarrow\cP^n\quad,\quad
\mu e_i=P^i\chi P^{n-i}
\ee
Theorem \ref{th-bm} implies that if $l=l(n)$ the elements in (\ref{basis})
lie in $\ker\mu$, so Propositions \ref{prop-basis} and \ref{prop-basis2}
imply that (\ref{basis}) with $l=l(n)$ is a basis for $\ker\mu$ and that
(\ref{basis2}) is a basis for $V_n\subseteq\cP^n$. Thus:

\begin{theorem} Any $l+1$ consecutive relations from the set 
{\rm(\ref{bm})} with $l=l(n)$ form a basis of relations among the elements of
$\{P^i\chi P^{n-i}:0\leq i\leq n\}$. The set 
$\{P^i\chi P^{n-i}:l(n)+1\leq i\leq n\}$ is a basis for $V_n$.
\label{main-theorem}
\end{theorem}

\noindent
{\bf Acknowledgement.} We thank Richard Stanley for the slick proof of 
Proposition \ref{prop-basis}.

\section{Independence of the relations}

We wish to show that (\ref{basis}) is a linearly independent set.
Regard elements of $\FF_p^{n+1}$ as column vectors, and arrange
the $l+1$ vectors in (\ref{basis}) as columns in a matrix, 
which we claim is of rank 
$l+1$. The top square portion is the mod $p$ reduction of the 
$(l+1)\times(l+1)$ integral Toeplitz matrix $A_l(m)$ with $(i,j)$th entry 
\[
\binom{m+j-i}{l}\quad,\quad 0\leq i,j\leq l
\]

\begin{lemma} $\det A_l(m)=1$.
\end{lemma}

\noindent
{\em Proof.} By induction on $m$. Since $\binom{-1}{l}=(-1)^l$ and
$\binom{-1+j}{l}=0$ for $0<j\leq l$, $A_l(-1)$ is lower triangular
with determinant $((-1)^l)^{l+1}=1$. Now we note the identity
\[
BA_l(m)=A_l(m+1)
\]
where 
\[
B=\left[\begin{array}{ccccc}
\binom{l+1}{1} & -\binom{l+1}{2} & \cdots & (-1)^{l-1}\binom{l+1}{l} & 
(-1)^l\binom{l+1}{l+1} \\
1 & 0 & \cdots & 0 & 0 \\
0 & 1 & \cdots & 0 & 0\\
\vdots & \vdots & & \vdots & \vdots \\
0 & 0 & \cdots & 1 & 0 
\end{array}\right].
\]
The matrix identity is an expression of the binomial identity
\be
\sum_k(-1)^k\binom{l+1}{k}\binom{n-k}{l}=0
\label{bin1}
\ee
(taking $n=m+1-j$ and $k=j+1$). Since $\det B=1$, the result follows 
for all $m\in\ZZ$. $\Box$

For completeness, we note that (\ref{bin1}) is the case $m=l+1$ of:

\begin{lemma}
$\displaystyle{\sum_k(-1)^k\binom{m}{k}\binom{n-k}{l}=\binom{n-m}{l-m}}$.
\end{lemma}

\noindent
{\em Proof.} The defining identity for binomial coefficients implies the case
$m=1$, and also that both sides satisfy the recursion 
$C(l,m,n)-C(l,m,n-1)=C(l,m+1,n)$. $\Box$

\section{Independence of the operations}

We will prove Proposition \ref{prop-basis2} by studying how 
$P^i\chi P^{n-i}$ pairs against elements in $\mathcal P_*$, the dual of the 
Hopf algebra of Steenrod reduced powers. According to Milnor \cite{milnor}, 
with our grading conventions 
\[
\cP_*=\FF_p[\xi_1,\xi_2,\ldots] \quad,\quad
|\xi_j|=\frac{p^j-1}{p-1}\,,
\]
and
\be
\Delta\xi_k=\sum_{i+j=k}\xi_i^{2^j}\tensor\xi_j\,.
\label{milnor-diagonal}
\ee
For a finitely nonzero sequence of nonnegative integers
$R=(r_1,r_2,\ldots)$ write $\xi^R=\xi_1^{r_1}\xi_2^{r_2}\cdots$ and let
$\lVert R\rVert=r_1+pr_2+p^2r_3+\cdots$ and 
\[
|R|=|\xi^R|=
r_1+\left(\frac{p^2-1}{p-1}\right)r_2+\left(\frac{p^3-1}{p-1}\right) 
r_3+\cdots\,.
\]
The following clearly implies Proposition \ref{prop-basis2}.

\begin{proposition} For any integer $n>0$ there exist sequences
$R_{n,j}$, $0\leq j\leq n-l(n)-1$, such that $|R_{n,j}|=n$ and
\[
\langle P^i\chi P^{n-i},\xi^{R_{n,j}}\rangle=\left\{
\begin{array}{cll}
\pm1 & \hbox{for} & i=n-j \\
0 & \hbox{for} & i>n-j\,. 
\end{array}\right.
\]
\label{basis3}
\end{proposition}

\noindent
The starting point in proving this is the following result of Milnor
\cite{milnor}.

\begin{lemma} [\cite{milnor}, Lemma 10] 
$\langle\chi P^n,\xi^R\rangle\neq0$ for all sequences $R$ with $|R|=n$. 
\end{lemma}

In the basis of $\cP$ dual to the monomial basis of $\cP_*$, the element
corresponding to $\xi_1^i$ is $P^i$. Since the diagonal in $\cP_*$ is dual
to the product in $\cP$, it follows from (\ref{milnor-diagonal}) that 
\[
\langle P^i\chi P^{n-i},\xi^R\rangle=\left\{
\begin{array}{cll}
\pm1 & \hbox{for} & i=\lVert R\rVert\\
0 & \hbox{for} & i>\lVert R\rVert\,.
\end{array}
\right.
\]

So we wish to construct sequences $R_{n,j}$, for $l(n)+1\leq j\leq n$,
such that $|R_{n,j}|=n$ and $\lVert R_{n,j}\rVert=j$. We deal first with the
case $j=l(n)+1$.

\begin{proposition} For any $n\geq0$ there is a sequence 
$M=(m_1,m_2,\ldots)$ such that \\
{\rm(1)} $|M|=n$, \\
{\rm(2)} $0\leq m_i\leq p$ for all $i$, and \\
{\rm(3)} If $m_j=p$ then $m_i=0$ for all $i<j$. \\
For any such sequence, $\lVert M\rVert=l(n)+1$. 
\label{M}
\end{proposition}

\noindent
{\em Proof.} Give the set of sequences of dimension $n$ 
the right-lexicographic order. 
We claim that the maximal sequence satisfies the hypotheses. 

Suppose that $R=(r_1,r_2,\ldots)$ does not satisfy the hypotheses.
If $r_1>p$ then the sequence $(r_1-(p+1),r_2+1,r_3,\ldots)$ is larger. 
If $r_j>p$, with $j>1$, then the sequence
$(r_1,\ldots,r_{j-2},r_{j-1}+p,r_j-(p+1),r_{j+1}+1,r_{k+2},\ldots)$ is larger.
This proves (2). 
To prove (3), suppose that $r_j=p$ with $j>1$, and suppose that some
earlier entry is nonzero. Let $i=\min\{k:r_k>0\}$. If $i=1$, then the 
sequence $(r_1-1,r_2,\ldots,r_{j-1},0,r_{j+1}+1,r_{j+2},\ldots)$ is larger.
If $i>1$, then $S$ with $s_k=0$ for $k<i-1$ and $i\leq k\leq j$,
$s_{i-1}=p$, $s_{j+1}=r_{j+1}+1$, and $s_k=r_k$ for $k>j+1$, is larger.

Let $M$ be a sequence satisfying (1)--(3), and write $l=\lVert M\rVert-1$. 
To see that $l=l(n)$  we must show that
\be
p(l+1)-\alpha(l+1)\geq(p-1)n
\label{greater}
\ee
and
\be
pl-\alpha(l)<(p-1)n\,.
\label{less}
\ee
The {\em excess} $e(R)$ is the sum of the entries in $R$, so that 
$p\lVert R\rVert-e(R)=(p-1)|R|$.  The $p$-adic representation of a number
minimizes excess, so for any sequence $R$ we have 
$e(R)\geq\alpha(\lVert R\rVert)$ and hence 
$p\lVert R\rVert-\alpha(\lVert R\rVert)\geq(p-1)|R|$: so (\ref{greater})
holds for any sequence. 

To see that (\ref{less}) holds for $M$, let 
$j=\min\{i:m_i>0\}$, so that $(p-1)n=(p^j-1)m_j+(p^{j+1}-1)m_{j+1}+\cdots$ and
$l+1=p^{j-1}m_j+p^jm_{j+1}+\cdots$. The hypotheses imply that 
$l$ has $p$-adic expansion
\[
(1+\cdots+p^{j-2})(p-1)+p^{j-1}(m_j-1)+p^jm_{j+1}+\cdots,
\]
so 
\[
\alpha(l)=(j-1)(p-1)+(m_j-1)+m_{j+1}+\cdots
\]
from which we deduce
\[
pl-\alpha(l)=(p-1)(n-j)<(p-1)n\,.
\]

\begin{corollary}
The function $l(n)$ is weakly increasing.
\end{corollary}

\noindent
{\em Proof.}
Let $M$ be a sequence satisfying the conditions of Proposition \ref{M},
and note that the sequence $R=(1,0,0,\ldots)+M$ has $|R|=n+1$ and 
$\lVert R\rVert=\lVert M\rVert+1=l(n)+2$. If $p$ does not occur in
$M$, then $R$ satisfies the hypotheses of the proposition (in degree $n+1$)
and hence $l(n)+1\leq l(n+1)+1$.
If $p$ does occur in $M$, then the moves described above will
lead to a sequence $M'$ satisfying the hypotheses. None of the moves
decrease $\lVert-\rVert$, so $l(n)+1\leq l(n+1)+1$. $\Box$

\begin{remark} {\rm Properties (1)--(3) in fact determine $M$ uniquely.} 
\end{remark}

\noindent
{\em Proof of Proposition \ref{basis3}.} Define $R_{n,l(n)+1}$ to be a sequence
$M$ as in Proposition \ref{M}. Then inductively define
\[
R_{n,j}=(1,0,0,\ldots)+R_{n-1,j-1}\quad\hbox{for}\quad l(n)+1<j\leq n\,.
\]
This makes sense by monotonicity of $l(n)$, and the elements clearly
satisfy $|R_{n,j}|=n$ and $\lVert R_{n,j}\rVert=j$. This completes the
proof. $\Box$

\end{document}